\documentclass[11pt]{article}
\usepackage{amsmath}
\usepackage{amssymb}
\usepackage{graphicx}

\numberwithin{equation}{section}

\title{Volume Optimization, Normal Surfaces
and Thurston's Equation on Triangulated 3-Manifolds   }
\author{Feng Luo\thanks{Partially Supported by a NSF grant}
}
\begin{document}

\maketitle


\begin{abstract}
\noindent We propose a finite dimensional variational principle on
triangulated 3-manifolds so that its critical points are related
to solutions to Thurston's gluing equation and Haken's normal
surface equation. The action functional is the volume. This is a
generalization of an earlier program by Casson and Rivin for
compact 3-manifolds with torus boundary.  Combining the result in
this paper and the work of Futer-Gu\'eritaud, Segerman-Tillmann
and Luo-Tillmann, we obtain a new finite dimensional variational
formulation of the Poncare-conjecture. This provides a step toward
a new proof the Poincar\'e conjecture without using the Ricci
flow.
\end{abstract}

\tableofcontents

\section{Introduction}\label{s-1}


\subsection{ The statement of the main theorem}
Given a closed triangulated  3-manifold or pseudo 3-manifold,
there are several linear and algebraic systems of equations
associated to the triangulation. Beside the homology theory, the
most prominent ones are Haken's theory of normal surfaces
(\cite{Ha}), Thurston's algebraic gluing equations for
constructing hyperbolic metrics (\cite{Th}) using hyperbolic ideal
tetrahedra, and the notion of angle structures (\cite{La},
\cite{Ri1}). The normal surface theory gives a parametrization of
essential surfaces in the manifold and solutions of Thurston's
equation produce hyperbolic cone metrics. Thurston used a solution
to Thurston's gluing equation to produce a complete hyperbolic
metric on the figure-8 knot complement in the earlier stage of
formulating his geometrization conjecture. The notion of
(Euclidean) angle structures, introduced by Casson, Rivin and
Lackenby  for 3-manifolds with torus boundary, is a linearized
version of Thurston's equation.

The goal of this paper generalizes the notion of angle structures
introduced by Casson, Rivin \cite{Ri1} and Lackenby \cite{La} to
the circle-valued angle structure (or $\bold S^1$-angle structure
for short) on any closed triangulated pseudo 3-manifold $(M, \bold
T)$. Using the method introduced in \cite{LT}, we show that
circle-valued angle structures always exist on any $(M, \bold T)$.
 Furthermore, the space of all $\bold
S^1$-angle structures on $(M, \bold T)$, denoted by $SAS(M, \bold
T)$, is shown to be a closed smooth manifold  (proposition 2.6).
Each $\bold S^1$-angle structure has a natural volume given by the
Milnor-Lobachevsky function. This defines a continuous, but not
necessary smooth volume function $\bold V :  SAS(M, \bold T) \to
\bold R$.  Since the space $SAS(M, \bold T)$ is compact, the
volume function $\bold V$ has a maximum point.  Our main result is
the following.
\\
\\
{\bf Theorem 1.1} \it  Suppose $(M, \bold T)$ is a triangulated
closed orientable pseudo 3-manifold.  Let $p$ be a maximum point
of the volume function $\bold V :SAS(M, \bold T)  \to \bold R$.

(a) If $p$ is a smooth point for $\bold V$, then $p$ produces a
solution to the generalized Thurston gluing equation.

(b) If $p$ is a non-smooth point for $\bold V$, then $p$ produces
a solution to Haken's normal surface equation with exactly one or
two non-zero quadrilateral coordinates. \rm
\medskip

In the sequel, we call a solution to Haken's equation in part (b)
of theorem 1.1 a \it 2-quad-type solution. \rm   Recall that a
triangulation of $M$ is called \it minimal \rm if it has the
smallest number of tetrahedra among all triangulations of $M$.
 It is conceivable
that the existence of a 2-quad-type solution on a minimally
triangulated 3-manifold puts constrains on the topology of the
manifold.  This is indeed the case. In our recent joint paper with
S. Tillmann \cite{LT2}, using the work of Jaco-Rubinstein
\cite{JR} and Futer-Gu\'eritaud, we proved that

\medskip

\noindent {\bf Theorem 1.2} (Luo-Tillmann \cite{LT2}) \it Suppose
$(M, T)$ is a minimally triangulated orientable closed 3-manifold
so that the volume function $\bold V :  SAS(M, \bold T) \to \bold
R$ has a non-smooth maximum point. Then,

(a) $M$ is reducible, or

(b) $M$ is toroidal, or

(c) $M$ is a Seifert fibered space, or

(d) $M$ contains the connected sum $\#_{i=1}^3 RP^2$ of three
copies of the projective plane. \rm

\medskip
Theorems 1.1 and 1.2 prompt us to make
\\
\\
{\bf Conjecture 1} \it Suppose $(M, \bold T)$ is a minimally
triangulated closed irreducible orientable 3-manifold so that all
maximum points of $\bold V: SAS(\bold T) \to \bold R$ are smooth
for $\bold V$.  Then $(M, \bold T)$ supports a solution to
Thurston's gluing equation. \rm
\medskip

We thank Ben Burton and Henry Segerman for providing data which
help us to formulate conjecture 1.  As we will see in the next
subsection, conjecture 1 for simply connected manifold is
equivalent to the Poincar\'e conjecture.

A weaker version of conjecture 1 is conjecture 2 in \S5. It does
not involve a maximization process and deals only with solutions
to Thurston's equation and Haken's equation. It is shown in \S5
that for a simply connected 3-manifolds, conjecture 2 is
equivalent to the Poincar\'e conjecture.

\subsection{Thurston's gluing equation}

Recall that a closed pseudo 3-manifold is the quotient space of a
disjoint union of tetrahedra so that their codimension-1 faces are
identified in pairs by affine homeomorphisms. In particular, a
closed 3-manifold is a pseudo 3-manifold.  The \it generalized
Thurston gluing equation \rm  associated to a triangulated
oriented pseudo 3-manifold is defined as follows. Assign each edge
in each tetrahedron in the triangulation $\bold T$ a complex
number $z \in \bold C-\{0,1\}$. The assignment is said to satisfy
\it the generalized Thurston gluing equation \rm if

\medskip

\noindent (a) opposite edges of each tetrahedron have the same
assignment;

\noindent (b) the three complex numbers assigned to three pairs of
opposite edges in each tetrahedron are  $z$, $\frac{1}{1-z}$ and
$\frac{z-1}{z}$ subject to an orientation convention; and

\noindent (c) for each edge $e$ in the triangulation, if \{$z_1,
..., z_k$\} is the set of all complex numbers assigned to the edge
$e$ in various tetrahedra adjacent to $e$,
 then
\begin{equation}  \prod_{i=1}^k z_i = \pm 1.
\end{equation}

\begin{figure}
 \centering
\includegraphics[scale=0.5]{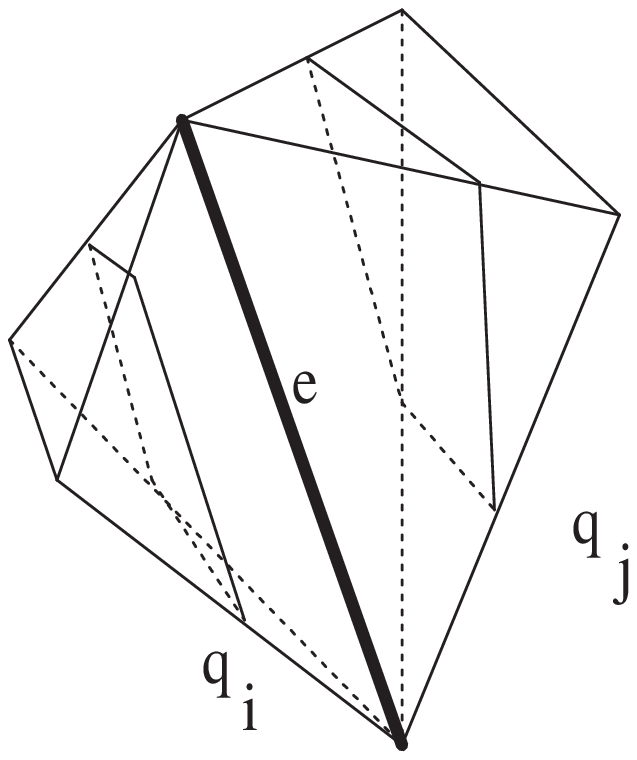}
\centerline{Figure 1.1}

\end{figure}

If the right-hand-side of (1.1) are 1 for all edges, we call the
assignment satisfying \it Thurston's gluing equation \rm (or \it
Thurston equation \rm for short).

We would like to emphasize that Thurston's equation and its
solutions are well defined on any closed oriented triangulated
pseudo 3-manifolds.  The most investigated cases are solutions of
Thurston equation on an ideal triangulated 3-manifold with torus
boundary so that the complex numbers $z$ are in the upper-half
plane  (see for instance \cite{Th}, \cite{Til2}, \cite{Dun},
\cite{PW} and many others). We intend to study Thurston's equation
in the most general setting. Even though a solution to Thurston
equation in the general setting does not necessary produce a
hyperbolic structure, one can still obtain some important
information from it.  For instance, it
 was observed in \cite{Yo} (see also \cite{NZ},
\cite{ST}) that each solution of Thurston's equation produces a
representation of the fundamental group of the 3-manifold with
vertices of triangulation removed to $PSL(2, \bold C)$.  A
simplified  version of a theorem of Segerman-Tillmann \cite{ST}
states that if $(M, \bold T)$ is a one-vertex triangulation of a
closed 3-manifold so that $\bold T$ supports a solution to
Thurston equation, then each edge in $\bold T$, considered a loop
in $M$, is homotopically essential in $M$. In particular, any
one-vertex triangulation of a simply connected 3-manifold  cannot
support a solution to Thurston equation.

Using this theorem of Segerman-Tillman and theorem 1.2, we can
deduce the Poincar\'e conjecture from conjecture 1 as follows.
Suppose $M$ is a simply connected closed 3-manifold. By the
Kneser-Milnor prime decomposition theorem, we may assume that $M$
is irreducible. Take a minimal triangulation $\bold T$ of $M$. By
the work of Jaco-Rubinstein on 0-efficient triangulation, we may
assume that $\bold T$ has only one vertex, i.e., each edge is a
loop. By Segerman-Tillmann's theorem above, we see that $(M, \bold
T)$ cannot support a solution to Thurston equation. By conjecture
1, there exists a non-smooth maximum volume $\bold S^1$-angle
structure.  By theorem 1.2, the minimality of $\bold T$ and
irreducibility of $M$, we conclude that $M = \bold S^3$.



Theorem 1.1 is a special case of theorem 3.2 in \S 3 where one
shows that under the same assumption as theorem 1.1 there are
either solutions to generalized Thurston's equation so that (a),
(b) hold, and for any edge $e$,
\begin{equation} \prod_{i=1}^k z_i = \pm k(e) \end{equation}
 where $k(e) \in \bold S^1$
is a given function satisfying (2.16) and (2.17) or there exists a
2-quad-type solution to Haken's equation. Theorem 1.1 provides
some evidence relating normal surface theory to representations.

A potential application of theorem 1.1 is to construct hyperbolic
metrics on closed 3-manifolds. Namely, if a maximum volume $\bold
S^1$-angle structure produces a solution to Thurston's equation so
that (i) the right-hand-side of Thurston's equations (1.1) are 1
and (ii) the maximum volume is the Gromov norm of  $M$ divided by
the volume of the ideal regular tetrahedron, then the associated
representation will likely produce a hyperbolic metric on $M$
(\cite{Dun}, \cite{Fra}). In our recent work with Tillmann and
Yang \cite{LTY}, we have shown that for a closed hyperbolic
3-manifold, there exists a triangulation and a solution to
Thurston's equation so that the above conditions (i) and (ii)
hold.

\medskip

\subsection{Remarks}

We remark that all results in the paper can be generalized without
difficulties to compact pseudo 3-manifolds with boundary. The
simplest way to treat them is by taking the doubling construction.
For simplicity, we will not state the corresponding theorems for
pseudo manifolds with boundary.

Using volume optimization on the space of angle structures to find
hyperbolic structures was carried out successfully by F.
Gu\'eritaud in \cite{Gu} for for punctured torus bundles over the
circle.  Our method is similar to that of \cite{Gu} in a different
setting.

The paper is organized as follows. In \S2, we revisit the theory
of normal surfaces and spun normal surfaces. In \S3, we recall
Thurston's work on gluing hyperbolic metrics and the volume of
angle structures. Theorem 1.1 is proved in \S4.  Some open
questions will be discussed in \S5.

We would like to thank S. Tillmann, D. Futer, F. Gueritaud, W.
Jaco, and Ruifeng Qiu for discussions.  The work is partially
supported by the NSF.

\section{The theory of normal surfaces revisited}

The  normal surface theory, developed by Haken in 1950's, is a
beautiful chapter in 3-manifold topology. In late 1970's, Thurston
introduced the notion of spun normal surfaces and used it to study
3-manifolds.  There are works by Tollefson, Kang-Rubinstein,
Tillmann,  Thurston, Jaco and others which characterize spun
normal surfaces using Haken's normal coordinates. It turns out a
spun normal surface is most conveniently described in terms of the
tangent vectors to $\bold S^1$-angle structures.  In fact, the two
system of linear equations for the tangential angle structures and
the spun normal surfaces are dual to each other. This observation,
which is implicit in the work of \cite{Tol}, \cite{KR}, and
\cite{Til}, is very useful for us in \S4 to relate the critical
points of the volume functional with the normal surfaces.

We will revisit the normal surface theory and follow the
expositions in \cite{JT} and \cite{Til} closely in this section.
Some of the notations used in the section are new.

\subsection{ Triangulations of closed pseudo 3-manifolds and
normal surfaces}

\medskip
Let $X$ be a union of finitely many disjoint oriented Euclidean
tetrahedra. The collection of all faces of tetrahedra in $X$ is a
simplicial complex $\bold T^*$ which is a triangulation of $X$.
Identify codimension-1 faces in $X$ in pairs  by affine
homeomorphisms. The quotient space $M$ is a closed pseudo
3-manifold with a triangulation $\bold T$ whose simplexes are the
quotients of simplexes in $\bold T^*$. See \cite{JR} for more
information.

Note that in this definition of triangulation, we do not assume
that simplexes in $\bold T$ are embedded in $M$. For instance, it
may well be that $\bold T$ has only one vertex. Furthermore, the
non-manifold points in $M$ are either the vertices or the centers
of the edges. If we require the affine homeomorphisms used in the
gluing for $M$ be orientation-reversing, then the pseudo manifold
$M$ is oriented and non-manifold points of $M$ are contained in
the vertex set  $V$. Let $N$ be the compact 3-manifold obtained
from $M$ with a small open regular neighborhood of $V$ removed.
Then we call $\{ \sigma \cap N | \sigma \in \bold T\}$ an \it
ideal triangulation \rm of $N$.

\medskip
According to Haken \cite{Ha}, a \it normal arc \rm in $X$ is an
embedded arc in a triangle face so that its end points are in
different edges and a \it normal disk \rm in $X$ is an embedded
disk in a tetrahedron so that its boundary consists of 3 or 4
normal arcs. These are called \it normal triangles \rm and \it
normal quadrilaterals \rm respectively.

\begin{figure}

\centering
\includegraphics[scale=0.6]{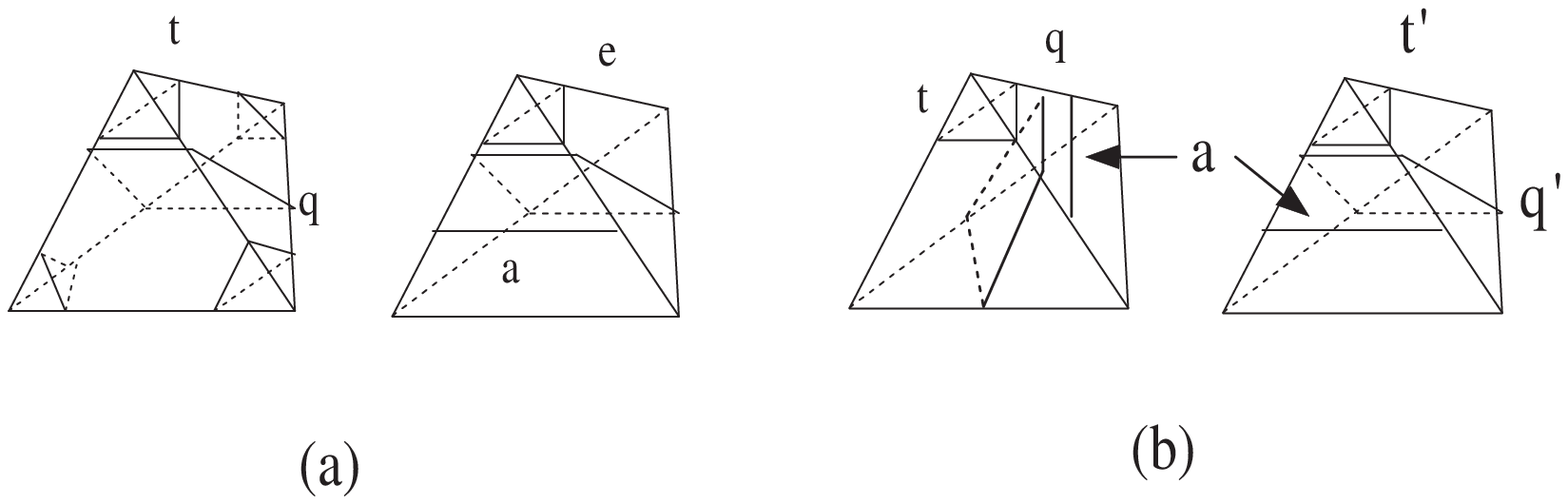}

\centerline{Figure 2.1}

 \end{figure}



The projections of normal arcs and normal disks from $X$ to $M$
constitute normal arcs and normal disks in the triangulated space
$(M, \bold T)$.  A surface $S$ in $M$  (or $X$) is called \it
normal \rm with respect to the triangulation $\bold T$ (or $\bold
T^*$) if for each tetrahedron $\sigma$ in the triangulation, the
intersection $S \cap \sigma$ consists of normal disks.
 A \it normal isotopy \rm is an isotopy of the
ambient space $X$ or $M$ which leaves each simplex invariant.
Normal arcs and disks will be considered up to normal isotopy. In
each tetrahedron, there are four normal triangles and three normal
quadrilaterals up to normal isotopy. We use $\triangle$, $\Box$
and $\bold A$ to denote the sets of all normal isotopy classes of
normal triangles, quadrilaterals and normal arcs in the
triangulation $\bold T$.  Since the set of all normal isotopy
classes of normal quadrilaterals (and normal triangles) in $\bold
T^*$ is the same as $\Box$ (and $\triangle$). We will also use
$\Box$ and $\triangle$ to denote the sets of all normal isotopy
classes of normal quadrilaterals and normal triangles in  $\bold
T^*$. In the sequel, we will use both "triangle" and
"quadrilateral" for the normal isotopy class of a triangle and a
quadrilateral.

Let $V, E, F, T$ be the sets of all vertices, edges, triangles and
tetrahedra in $\bold T$.  The set of all edges and tetrahedra in
$\bold T^*$ will be denoted by $E^*$ and $T^*$.  We consider $E$
as the set of equivalence classes of elements in $E^*$, i.e., $
E=\{ [x] | x \in E^*$ where $x$ and $y$ in $E^*$ are equivalent if
they are identified in $X$\}.

If $x, y \in V \cup E \cup F \cup T$, we use $ x > y$ to denote
that $y$ is a face of $x$. We use $|Y|$ to denote the cardinality
of a set $Y$.

There are relationships among various sets $V, E, F, T, \triangle,
\Box, \bold A$. These incidence relations, which will be recalled
below, are the basic ingredient for defining linear and algebraic
equations on $\bold T$.

Take $ t \in \triangle$, $a \in \bold A$, $q \in \Box$, $e \in E$,
and $\sigma \in T$.  We use $ a< t$ (and $ a < q$) if there exist
representatives $x \in a$, $y \in t$ (and $z \in q$) so that $x$
is an edge of $y$ (and $z$). We use $t \subset \sigma$ and $q
\subset \sigma$ to denote that representatives of $t$ and $q$ are
in the tetrahedron $\sigma$. In this case, we say the tetrahedron
$\sigma$ contains $t$ and $q$.

The index $i(q, e)$ is an integer $0,1$ or $2$ defined as follows.
Given $q \in \Box$ and $e^* \in E^*$, let $i(e^*, q)$ be 1 if $q,
e^*$ lie in the same tetrahedron so that $ q \cap e^* =\emptyset$
and let $i(e^*, q) =0$ in all other cases. If $e \in E$ and $q \in
\Box$, then the  index $i(e,q) = \sum_{ e^* \in e} i(e^*, q)$.

The index $i(t, e)$ is the number of edges $e^* \in e$ so that $t$
has a vertex in $e^*$.  See figure 2.2.

\medskip

\begin{figure}

\centering
\includegraphics[scale=0.5]{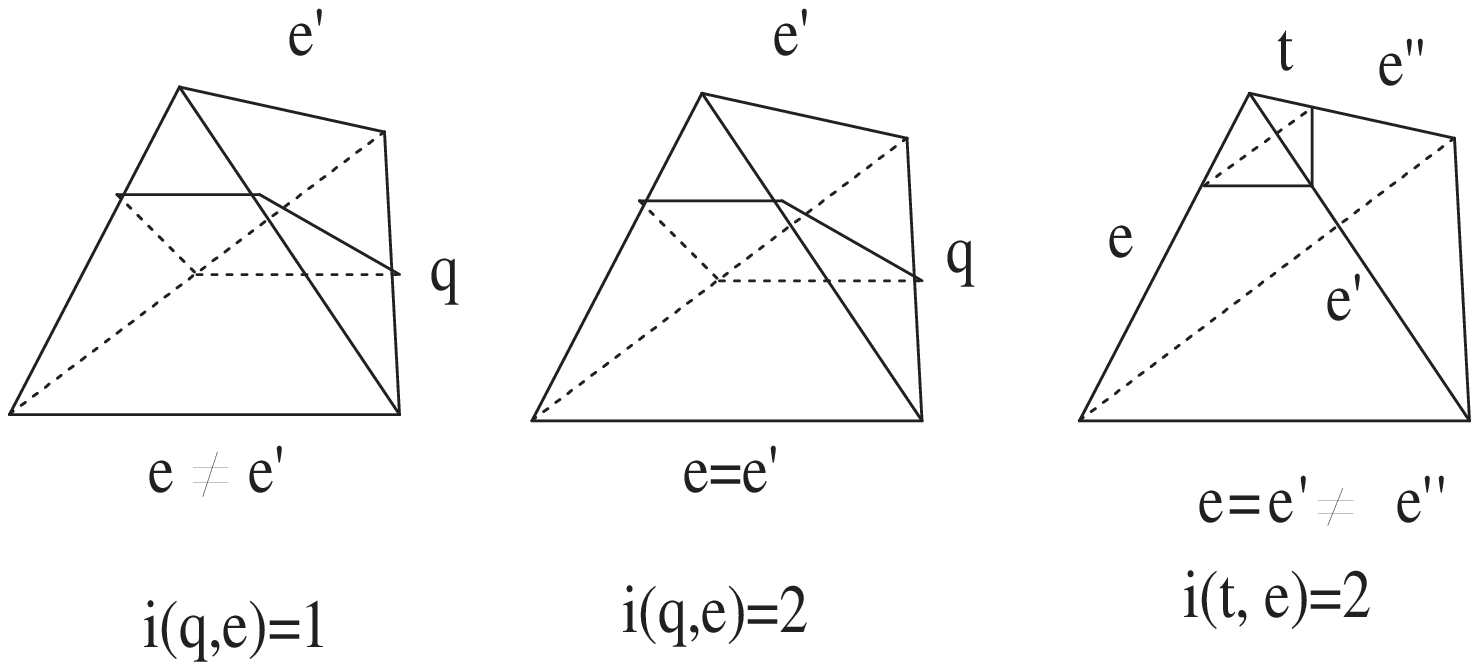}
\centerline{Figure 2.2}

 \end{figure}

We remark that if $\bold T$ is a simplicial triangulation, then
the indices $i(t, e)$ and $i(q, e)$ take only values $0, 1$.

As a convention, in the sequel, we will always use $\sigma$, $e$
and $q$ to denote a tetrahedron, an edge and a quadrilateral in
the triangulation $\bold T$ respectively.

\subsection{ Normal surface equation and Kang-Rubinstein
basis}

The normal surface equation is a system of linear equations
defined in the space $\bold R^{\triangle} \times \bold R^{\Box}$
introduced by W. Haken \cite{Ha}.  It is defined as follows. For
each normal arc $ a \in \bold A$, suppose $\sigma, \sigma'$ are
the two tetrahedra adjacent to the triangular face which contains
$a$. Then there is a homogeneous linear equation for  $ x \in
\bold R^{\triangle}\times \bold R^{\Box}$ associated to $a$:

\begin{equation} \label{2.1}
 x (t) + x(q) = x(q') + x(t')
 \end{equation}
 where  $t, q \subset \sigma$, $t', q' \subset \sigma'$ and $t, t',
 q, q' > a$. See figure 2.1(b). Let $\bold S_{ns}$ be the space of all solutions to
 (2.1) as $a$ runs over all normal arcs.

A basis of the solution space $\bold S_{ns}$ to equations (2.1)
was found by Kang-Rubinstein \cite{KR}. To state it, let us
introduce
 one more notation.  Given a finite set $Z$, the \it standard basis
 \rm of the vector space $\bold R^Z$ will be denoted by $\{ z^* |
 z \in Z\}$ where $ z^*(z) =1$ and $z^*(z') =0$ if $z' \in Z-\{z\}
 $. We  give $\bold R^Z$ the inner product so that $\{z^* | z
 \}$ forms an orthonormal basis.
Now for each $\sigma \in T$ and $e \in E$, define the vectors
$W_{\sigma}, W_e \in \bold R^{\triangle} \times \bold R^{\square}$
as follows,
\begin{equation}
W_{\sigma} = \sum_{q \in \Box, q \subset \sigma} q^*   -\sum_{ t
\in \triangle, t \subset \sigma} t^*\end{equation}

and \begin{equation} W_e =\sum_{ q \in \Box} i(q, e) q^* -\sum_{ t
\in \triangle} i(t,e) t^*.
\end{equation}

A basic theorem proved in \cite{KR} says,

\bigskip
\noindent {\bf Theorem 2.1}(Kang-Rubinstein). \it For any
triangulated closed pseudo 3-manifold, the set $\{W_x | x \in E
\cup T\}$ forms a basis of the solution space $\bold S_{ns}$ of
the normal surface equation. \rm

\bigskip
For the convenience of the reader, an alternative interpretation
of Kang-Rubinstein's proof is given in the appendix.f

\bigskip

\subsection{Spun normal surfaces and tangential angle structure}

Given  $x \in \bold R^{\triangle} \times \bold R^{\Box}$,  we will
call $x(t)$ and $x(q)$ the $t$-coordinate and $q$-coordinate
(triangle and quadrilateral coordinates) of $x$. Spun normal
surface theory addresses the following question, first
investigated by Thurston \cite{Th}. Given a vector $z \in \bold
R^{\Box}$, when does there exist a solution $x \in \bold S_{ns}$
to (2.1) whose projection to $\bold R^{\square}$ is $z$?
Geometrically, it asks if a given finite set of normal
quadrilaterals can be realized as the set of all normal
quadrilaterals in a normal surface. The question was completely
solved in the work of \cite{Tol}, \cite{KR}, \cite{Til} and
\cite{Ja2}. The results of Kang-Rubinstein and Tillmann are more
general and give solutions to the projections of not necessary
closed normal surfaces.

The purpose of this section is to interpret a weak version of
their work in terms of tangential angle structures.

\bigskip
\noindent {\bf Definition 2.1} A \it tangential angle structure
\rm on a triangulated pseudo 3-manifold $(M, \bold T)$ is a vector
$ x \in \bold R^{\Box}$ so that

 for each tetrahedron $\sigma \in T$, \begin{equation} \sum_{ q
\in \Box, q \subset \sigma} x(q)=0, \end{equation}

and for each edge $e \in E$, \begin{equation} \sum_{q \in \Box}
i(q, e) x(q) =0. \end{equation}

The linear space of all tangential angle structures on $(M, \bold
T)$ is denoted by $TAS(M, \bold T)$ or simply $TAS(\bold T)$.
\medskip

Recall that a (Euclidean type) angle structure, introduced by
Casson, Rivin and Lackenby, is a vector $ x \in \bold
R_{>0}^{\Box}$ so that for each tetrahedron $\sigma \in T$,
\begin{equation} \sum_{ q \in \Box, q \subset \sigma}
x(q)=\pi, \end{equation} and
 for each $e \in E$, \begin{equation} \sum_{q \in \Box} i(q, e)
x(q) =2\pi. \end{equation}  Thus one sees easily that a tangential
angle structure is a tangent vector to the space of all angle
structures.  In \cite{LT}, a \it generalized angle structure \rm
on $(M, \bold T)$ is defined as a vector $ x \in \bold R^{\Box}$
so that (2.6) and (2.7) hold. It is proved in \cite{LT} that a
generalized angle structure exists if and only if the euler
characteristic of each link $lk(v)$ is zero for $v \in V$.  We
will consider in this paper those $x \in \bold R^{\square}$ so
that the right-hand-side of (2.6) is in $\pi + 2\pi \bold Z$ and
the right-hand-side of (2.7) is in $2\pi \bold Z$. These will be
called \it $\bold S^1$-valued angle structures \rm on $\bold T$
and will be shown to exist on any closed pseudo 3-manifold using
the method introduced in \cite{LT}. Evidently $TAS(\bold T)$ is
the tangent space to $SAS(\bold T)$.

The following is a result proved by Tollefson (for closed
3-manifold case), Kang-Rubinstein and Tillmann for all cases. The
result was also known to Jaco \cite{Ja2}.

\medskip
\noindent {\bf Theorem 2.2} (\cite{Tol}, \cite{KR}, \cite{Til})
\it For a triangulated closed pseudo 3-manifold $(M, \bold T)$,
let $Proj_{\Box}: \bold R^{\triangle} \times \bold R^{\Box} \to
\bold R^{\Box}$ be the projection. Then
\begin{equation}
Proj_{\Box}(\bold S_{ns}) =TAS(\bold T)^{\perp}
\end{equation}
where $\bold R^{\Box}$ has the standard inner product so that
\{$q^* | q \in \Box\}$ is an orthonormal basis. \rm

\medskip
We remark that theorem 2.2 is not stated in this form in the work
of \cite{Tol}, \cite{KR}, \cite{Til}. This interpretation is due
to us.

\medskip
\noindent {\bf Proof.} Suppose $\bold R^n$ and $\bold R^m$ are
Euclidean spaces with the standard inner product and $A: \bold R^n
\to \bold R^m$ is a linear transformation with transpose $A^t:
\bold R^m \to \bold R^n$. Then it is well known that
$Im(A)=ker(A^t)^{\perp}$.  Define a linear map \begin{equation}
\label{ 2.9} A: \bold R^E \times \bold R^T \to \bold R^{\Box}
\end{equation}
by \begin{equation} A(h) = Proj_{\Box} ( \sum_{e \in E} h(e) W_e +
\sum_{\sigma \in T} h(\sigma) W_{\sigma}). \end{equation} By
definition, $Proj_{\Box}(\bold S_{ns}) = Im(A)$.

By (2.2) and (2.3),  we have
\begin{equation}
A(h)(q) =\sum_{ \sigma \in T, q \subset \sigma} h(\sigma) +
\sum_{e \in E} i(q,e) h(e).  \end{equation}

To understand tangential angle structures, we define a linear map
$ B: \bold R^{\Box} \to \bold R^E \times \bold R^T $ so that
\begin{equation}
B(x)(e) =\sum_{q \in \Box} i(q, e) x(q) \end{equation} and
\begin{equation}
B(x)(\sigma) =\sum_{ q \in \Box, q \subset \sigma} x(q).
\end{equation}

By definition, we have $TAS(\bold T)  = ker(B)$. We claim that $B
= A^t$, i.e., $(B(x), h) = (x, A(h))$ for all $x \in \bold
R^{\Box}, h \in \bold R^E \times \bold R^T$ where $(,) $ is the
standard inner product.

Indeed, by definition, we have

$$(B(x), h) = \sum_{e \in E} h(e)B(x)(e) + \sum_{ \sigma \in T}
h(\sigma)B(x)(\sigma)$$ $$ =\sum_{e \in E, q \in \Box} i(e,
q)x(q)h(e) + \sum_{\sigma \in T, q \in \Box, q \subset \sigma}
h(\sigma) x(q)$$
$$=\sum_{q \in \Box} x(q) \sum_{e \in E} i(q,e) h(e) + \sum_{ q
\in \Box} x(q) \sum_{ \sigma \in T, q \subset \sigma} h(\sigma)$$
$$=(x, A(h)).$$

Therefore, $TAS(\bold T)^{\perp} = ker(B)^{\perp} =Im(A) =
Proj_{\Box}(\bold S_{ns})$.  This ends the proof.
\bigskip

\noindent {\bf Corollary 2.3} (Tillmann [Ti]). \it (a)
$\dim(TAS(\bold T)) = |V| -|E| + 2|T| = \chi(M) + |T|$.

(b) $\dim(Proj_{\Box}(\bold S_{ns})) = -\chi(M) + 2|T|.$  \rm

\bigskip

\noindent {\bf Definition 2.2} Suppose $(M, \bold T)$ is a
triangulated closed pseudo 3-manifold. We say the triangulation
$\bold T$ is \it angle rigid \rm if there is $q \in \Box$ so that
$x(q)=0$ for all $x \in TAS(\bold T)$.  We say $\bold T$ is \it
2-angle rigid \rm if there exists a non-zero vector $(c_1, c_2)
\in \bold R^2$ and $q_1 \neq  q_2 \in \Box$ so that $c_1 x(q_1) +
c_2 x(q_2) =0$ for all $x \in TAS(\bold T)$.

\medskip
By definition, if $\bold T$ is angle rigid, then $x(q)$ is a
constant for all $\bold S^1$-angle structures $x$, i.e., the angle
at $q$ cannot be deformed.  If the triangulation $\bold T$ has an
edge $e$ of degree 1, then $\bold T$ is angle rigid at the
quadrilateral $q$ so that $i(q,e) \neq 0$. If $\bold T$ has an
edge $e$ of degree 2, then $\bold T$ is 2-angle rigid at the
quadrilaterals $q_1$ and $q_2$ so that $i(q_j, e) \neq 0$ for
$j=1,2.$

One simple consequence of Theorem 2.2 is,

\medskip
\noindent {\bf Corollary 2.4} \it Under the same assumption as in
theorem 2.2,

(a)  $(M, \bold T)$ is angle-rigid if and only if there exists an
embedded normal surface $\Sigma$ in $\bold T$ so that the surface
has exactly one normal quadrilateral type.

(b) $(M, T)$ is 2-angle rigid if and only if there exists a vector
 $ v \in \bold S_{ns} \cap (\bold Z
 ^{\triangle} \times \bold Z^{\Box})$ so that $Proj_{\Box}(v)$ is non-zero and has at most
 two non-zero coordinates.
 \rm

 \medskip
 To see part (a), if there exists a normal surface containing only one
 quadrilateral type
 $q \in \Box$, then its normal
 coordinate $x \in \bold R^{\triangle} \times \bold R^{\Box}$ is a
 vector so that $ Proj_{\Box}(x)  = k q^*\in TAS(\bold T)^{\perp}$ for
 some non-zero scalar $k$ and some $q \in \Box$.
 Thus $z(q) =0$ for all $z \in TAS(\bold T)$. Conversely, if there
 exists $q \in \Box$ so that $z(q) =0$ for all $z \in TAS(\bold
 T)$, then  $q^* \in TAS(\bold T)^{\perp}$. By theorem 2.2, $q^* =Proj_{\Box}(v)$
 for some $v \in \bold S_{ns}$. We
 may choose $v \in \bold Q^{ \triangle} \times \bold Q^{\Box}$
 since the linear equations (2.1) have integer coefficient and $q^*$ has integer coordinates. It
 follows some integer multiple $kv$ has non-negative integer
 q-coordinates. Now add to the vector $kv$ a positive integer multiples of the
 normal coordinates of the normal surfaces
 $lk(v)$, the link of the vertex $v \in V$, so the resulting vector has positive t-coordinates.
  We obtain a vector $u \in \bold S_{ns}
  \cap ( \bold Z_{\geq 0}^{\triangle} \times Z_{\geq 0}^{\Box})$ with exactly one
 non-zero q-coordinate. By the work of Haken, this vector $u$ is
  the normal coordinate of an embedded normal surface in $(M,
  \bold T)$.  The proof of (b) is similar and will be omitted.
  However, we are not able to conclude that $ v \in \bold Z_{\geq
  0}^{\triangle} \times \bold Z_{ \geq 0}^{\Box}$ in this case.

\medskip
\noindent {\bf Question.} Suppose $M$ is a non-compact complete
finite volume hyperbolic manifold and $\bold T$ is an ideal
geometric triangulation of $M$.
 Is $\bold T$ 2-angle rigid?

\medskip
\noindent
\subsection{Existence of $\bold S^1$-valued angle structures}

We begin with a definition which was also known to D. Futer and F.
Gueritaud.

\medskip

\noindent {\bf Definition 2.2}  Suppose $(M, \bold T)$ is a
triangulated closed pseudo 3-manifold. Let $k: E \to \bold S^1$ be
given. An \it $\bold S^1$-valued angle structure with curvature
$k$ \rm on $\bold T$ is a function $x: \square \to \bold S^1$ so
that for each tetrahedron $\sigma \in \bold T$,
\begin{equation}
\prod_{ q \in \square, q \subset \sigma}  x(q) =-1
\end{equation}
and for each edge $e \in E$,
\begin{equation}
\prod_{ q \in \square} x(q)^{i(q, e)} =k(e).
\end{equation}

The set of all $x \in (\bold S^1)^{\square}$ satisfying (2.14) and
(2.15) will be denoted by $SAS(\bold T, k)$.  The case that
$k(e)=1$ for all $e \in E$ is the most interesting one.  We use
$SAS(\bold T)$ to denote $SAS(\bold T, 1)$ where $1(e) =1$ for all
$e \in E$.

\medskip

For a complex number $w \in \bold C$, we use $\arg(w) \in [0,
2\pi)$ to denote its argument. If $x \in SAS(\bold T, k)$, by
taking $\arg(x(q))$, we can interpret an $\bold S^1$-valued angle
structure $x$ as a map from $\Box \to \bold R$  satisfying (2.6)
and (2.7) so that the right-hand side of (2.6) is in $2\pi \bold Z
+ \pi$ and the right-hand-side of (2.7) is in $2\pi \bold Z +
\arg(k(e))$.

\medskip
\noindent {\bf Lemma 2.5} \it If $SAS(\bold T, k) \neq \emptyset$,
then the function $k:E \to \bold S^1$ satisfies,

\begin{equation} \prod_{e \in E} k(e) =1 \end{equation}

and for each vertex $v \in V$,
\begin{equation}  \prod_{ e > v} k(e) =1. \end{equation} \rm

\medskip
Indeed, to see (2.16), using (2.14) and (2.15), we can write the
left-hand side of (2.16) as
 $$ \prod_{ q \in \square,  e \in E}  x(q)^{ i(q, e)}
 = \prod_{ \sigma \in T} \prod_{ q \in \sigma, e  < \sigma}
 x(q)^{ i(q, e)}
  = \prod_{ \sigma \in T}  (\prod_{ q \subset \sigma} x(q))^2
 =1$$
 due to $\sum_{ e \in E} i(q,e) =2$ for each $q$.

 To see (2.17), using (2.14), we can write the left-hand-side of (2.17)
 as,
 $$ \prod_{ e > v} \prod_{q \in \Box} x(q)^{i(q,e)}
 = \prod_{ \sigma \in T, \sigma > v}  (\prod_{ q \subset \sigma, e < \sigma,  e > v}
 x(q)^{i(q,e)}) =(-1)^N$$
 where $N$ is the number of normal triangles at the vertex $v$.
 This number $N$ is the same as the number of triangles in the
 link $lk(v)$. Since $lk(v)$ is a closed triangulated surface, $N$ is an
 even number. Thus (2.17) follows.

\medskip
 One can define the
similar notion of $\bold S^1$-valued angle structure on a closed
triangulated surface by assigning each angle of a triangle a
complex number of norm 1 so that the product of the complex
numbers in each triangle is $-1$. The \it curvature \rm at a
vertex is the product of all complex numbers assigned to the
angles at the vertex. For instance, if $(M^3, \bold T)$ is a
triangulated pseudo 3-manifold with an $\bold S^1$-valued angle
structure, then the vertex link $lk(v)$ has the induced $\bold
S^1$-angle structure. The identity (2.17) says that the product of
its curvatures at all vertices is 1.

The main result in this section says that (2.16) and (2.17) are
also sufficient.  This generalizes our earlier work with Tillmann
on 3-manifolds with torus boundary \cite{LT}. The method of the
proof of the proposition below is that of \cite{LT}.

\medskip
\noindent {\bf Proposition 2.6} (See also \cite{LT}). \it Given
any triangulated closed pseudo 3-manifold $(M, \bold T)$ and $k: E
\to \bold S^1$ satisfying (2.16) and (2.17), then $SAS(\bold T, k)
\neq \emptyset$. Furthermore, $SAS(\bold T, k)$ is a smooth closed
manifold of dimension $| \chi(M)| + |T|$. \rm

\medskip
\noindent {\bf Proof.} We may assume without loss of generality
that $M$ is connected.  Consider the Lie group homomorphism $F:
(\bold S^1)^{\square} \to (\bold S^1)^E \times (\bold S^1)^T$
given by
$$  F(z)(e) = \prod_{ q \in \square} z(q)^{i(q, e)} $$
and
$$ F(z)(\sigma) = \prod_{ q \in \square, q \subset \sigma}  z(q)$$
where $ z \in (\bold S^1)^{\square}$, $e \in E$ and $\sigma \in
T$.  The goal is to show that the point $t : E \cup F \to \bold
S^1$ given by $ t(e) =k(e)$ for $e \in E$ and $t(\sigma) =-1$ for
$\sigma \in T$ is in the image of $F$.

Suppose otherwise, that $t$ is not in the image of $F$. Since $F$
is a continuous group homomorphism from a torus to a torus, the
image of $F$ is a connected closed subgroup of $(\bold S^1)^E
\times (\bold S^1)^T$ which misses $t$. Thus there exists a
continuous group homomorphism $h: (\bold S^1)^E \times (\bold
S^1)^T \to \bold S^1$ so that $h(t) \neq 1$ and $h F$ is the
trivial homomorphism.

Each homomorphism from $(\bold S^1)^n$ to $\bold S^1$ is given by
a vector $(m_1, ..., m_n) \in \bold Z^n$, i.e., the homomorphism
sends $(x_1, ..., x_n) \in (\bold S^1)^n$ to $x_1^{m_1} ...
x_n^{m_n}$. Thus for the homomorphism $h$, there exists $\phi \in
\bold Z^E \times \bold Z^T$ so that for all $x \in (\bold S^1)^E
\times (\bold S^1)^T$,
$$ h(x) = \prod_{ e \in E} x(e)^{\phi(e)} \prod_{ \sigma \in T} x(\sigma)^{\phi(\sigma)}.$$

By the choice of $t$, we have $h(t) = \prod_{\sigma \in T}
(-1)^{\phi(\sigma)} \prod_{e \in E} k(e)^{\phi(e)}$. Thus $h(t)
\neq 1$ says that \begin{equation} \prod_{e \in E} k(e)^{\phi(e)}
\neq (-1)^{\sum_{ \sigma \in T} \phi(\sigma)}.
\end{equation} On the other hand, we will show that $\phi F $
being trivial implies that (2.18) is an equality. The
contradiction establishes the proposition.

Since the composition $hF$ is trivial,  for any $z \in (\bold
S^1)^{\square}$,

$$ h(F(z)) = (\prod_{e \in E} \prod_{q} z(q)^{i(q,e) \phi(e)})(
\prod_{\sigma \in T} \prod_{q \subset \sigma} z(q)^{\phi(\sigma)})
$$
$$= \prod_{ q \in \square} z(q)^{\sum_{\sigma, q \subset \sigma} \phi(\sigma)+ \sum_{ e \in E}
\phi(e) i(q, e)}.$$

By the assumption, $h(F(z))=1$ for all choice of $z \in (\bold
S^1)^{\square}$. Thus we obtain, for each $q \in \square$,
\begin{equation}  \sum_{\sigma, q \subset \sigma} \phi(\sigma) + \sum_{e} i(q, e)
\phi(e) =0. \end{equation}

Fix a tetrahedron $\sigma \in \bold T$, the above equation says
that the sum of $\phi(e) + \phi(e')$ of the values of $\phi$ at
two opposite edges $e, e'$ in $\sigma$ is independent of the
choice of $e, e'$. We will need to use the following lemma.

\medskip
\noindent
 {\bf Lemma 2.7} \it Suppose $a_{ij}=a_{ji} \in \bold Z$
where $ i \neq j \in \{1,2,3,4\}$ are six numbers so that
$$ a_{ij}+a_{kl}=c $$ is a constant independent of choice of indices $i,j,k,l$
where $\{i,j,k,l\}=\{1,2,3,4\}$. Then there exist $b_1, .., b_4
\in \bold Z/2 =\{ n/2 | n \in \bold Z\}$ so that
$$  a_{ij} = b_i + b_j$$ for all $i \neq j \in \{1,2,3,4\}$. \rm

\medskip
Indeed, $b_i = \frac{ a_{ij}+ a_{ik} - a_{jk}}{2}$ is independent
of the choices of $\{i,j,k\}$, $\{i, j, l\}$, or $\{i,k,l\}$ due
to the assumption on $a_{ij}+ a_{kl}=c$.

\medskip

\begin{figure}

\centering
\includegraphics[scale=0.5]{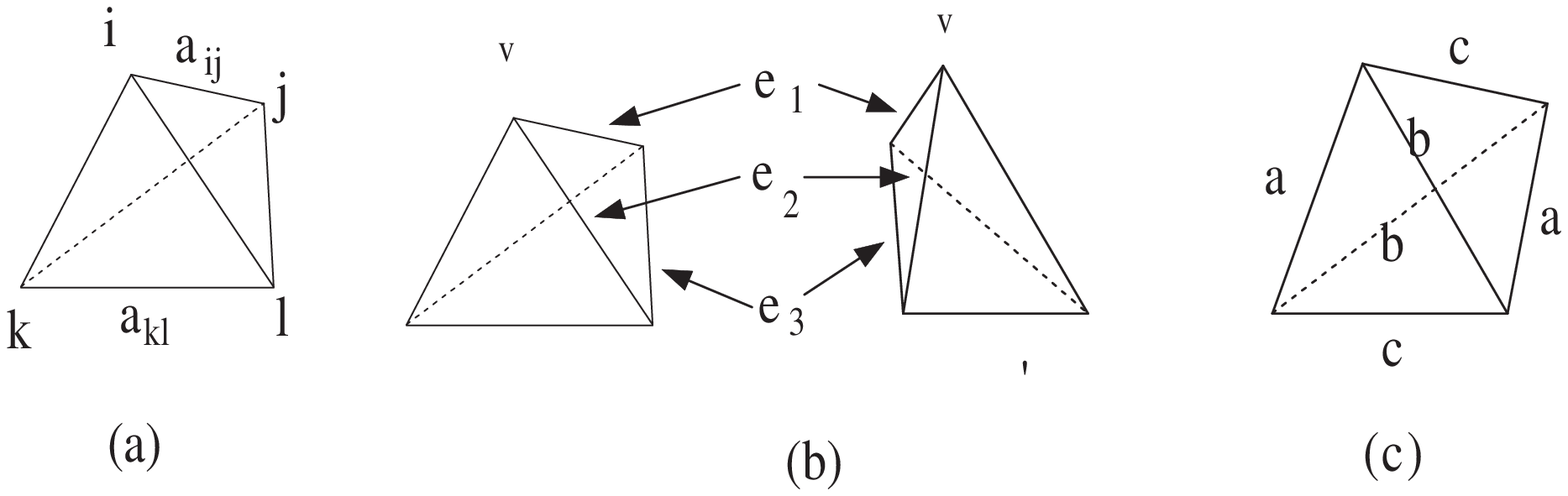}
\centerline{Figure 2.3: A topological interpretation of lemma 2.7
}

 \end{figure}

Thus, by the lemma, there exists a map $w: \{ \text {vertice of
$\sigma$} \} \to \bold Z/2$
 so that
\begin{equation} \phi(e) = w(v, \sigma) + w(v',
\sigma)\end{equation} where $v, v'$ are the end points of $e$. We
claim that $w(v, \sigma)$ is independent of the choice of
$\sigma$. Indeed, consider two tetrahedra $\sigma, \sigma'$
sharing a common triangular face $f$ (see figure 2.3(b)). Then for
three edges $e_1, e_2, e_3$ in $f$, we solve (2.20) and obtain
$$ w(v, \sigma) = w(v, \sigma') = \frac{ \phi(e_1) + \phi(e_2)
-\phi(e_3)}{2}$$ where $v$ is the vertex opposite to the edge
$e_3$ in $f$. It follows that $w(v, \sigma) =w(v, \sigma')$ is
independent of the choice of tetrahedra $\sigma$ and $\sigma'$
since $(M, \bold T)$ is a pseudo 3-manifold. Let  $w: V \to \bold
Z/2$ be the map so that

\begin{equation} \phi(e) = w(v) + w(v')\end{equation} where $v, v'$ are vertices of $e$.
We claim that either all $w(v)$'s are integers, or all of $w(v)$
are half-integers (i.e.,
 $k+1/2$ for some $k \in \bold Z$). Indeed, since $\phi(e)$ is an
 integer, it follows from (2.21) that either both $w(v), w(v')$ are in $\bold
 Z$, or both are in $\bold Z/2-\bold Z$.  Since the manifold $M$ is
 connected, it follows that either $w(v) \in \bold Z$ for all $v$,
 or $w(v) \in \bold Z/2 -\bold Z$ for all $v$.

We now claim that the sum $\sum_{ \sigma \in T} \phi(\sigma)$ has
to be an even integer. Indeed, by (2.19) and (2.21), $\phi(\sigma)
= -\sum_{ v < \sigma} w(v)$. Thus $$ \sum_{\sigma \in T}
\phi(\sigma) = -\sum_{ v \in V} w(v) (\sum_{ \sigma > v} 1)$$
$$= -\sum_{v \in V} w(v) |\{ \text{triangles in the link
lk(v)}\}|$$

For any triangulation of a closed surface, the number of triangles
in the triangulation has to be even. Thus if all $w(v) \in \bold
Z$, $\sum_{ \sigma \in T} \phi(\sigma)$ is even.  In the other
case, all $w(v) \in \bold Z/2 -\bold Z$. Thus
$$\sum_{\sigma} \phi(\sigma) = -\sum_{v \in V} \frac{1}{2} |\{
\text{triangles in the link lk(v)} \}|  \quad mod (2)$$ $$
=-\frac{1}{2} |\{ \text{ normal triangles in T}\}| \quad mod (2)$$
Now each tetrahedron has 4 normal triangles, thus the total number
of normal triangles in $T$ is a divisible by 4. Thus implies again
that $\sum_{\sigma \in T} \phi(\sigma) $ is an even number.

This implies that the right-hand-side of (2.18) is 1. We claim
that the left-hand-side of (2.18) is also equal to 1.  There are
two cases to be considered. In the first case, all $w(v)$'s are in
$\bold Z$. Then the left-hand-side of (2.18) becomes
$$ \prod_{e \in E} k(e)^{\sum_{ v < e} w(v)}  =\prod_{ v \in V} (
\prod_{ e > v} k(e) )^{ w(v)}$$ which is 1 due to (2.17).

In the second case that $ w(v) = W(v) + 1/2$ where $W(v) \in \bold
Z$ for all $v \in V$. We have $\phi(e) =1+\sum_{ v < e} W(v)$.
Thus the left-hand-side of (2.18) becomes
$$ \prod_{e \in E} k(e)^{1+\sum_{ v < e} W(v)}  =(\prod_{ v \in V}
\prod_{ e > v} k(e) ^{ W(v)} ) \prod_{ e \in E} k(e)$$ which is
again 1 due to (2.16) and (2.17). This contradict shows that
$SAS(\bold T, k) \neq \emptyset$.

Finally, since $SAS(\bold T, k) =F^{-1}(t)$ where $F$ is a Lie
group homomorphism, one concludes that $SAS(\bold T, k)$ is a
closed smooth manifold.

\bigskip
\section{Thurston's algebraic gluing equation and volume}

Thurston's equation mentioned in the introduction can be
conveniently rephrased in terms of the normal quadrilaterals in
the triangulation. It is based on the fact that a pair of opposite
edges in a tetrahedron is the same as the normal isotopy class of
the quadrilateral which is disjoint from the given edges.  We will
rewrite Thurston's equation in terms of quadrilaterals in this
section. In order to do so, we first recall the Neumann-Zagier
anti-symmetric bilinear form on $\bold R^{\Box}$. This bilinear
form appeared in the important work of Neumann-Zagier \cite{NZ}.
We assume that $(M, \bold T)$ is an oriented closed pseudo
3-manifold in this section so that each tetrahedron in $\bold T$
has the induced orientation.

\medskip
\subsection{Neumann-Zagier  Poisson
structure}

If $\sigma$ is an oriented Euclidean tetrahedron with edges from
one vertex labelled by $a,b,c$ so that the opposite edges have the
same labelling $a,b,c$ (see figure 2.3(c)), then the cyclic order
of edges $a,b,c$ viewed from each vertex depends only on the
orientation of the tetrahedron, i.e., is independent of the choice
of the vertices.  Now each pair of opposite edges  in the
tetrahedron corresponds to a normal isotopy class of quadrilateral
$q$ in $\sigma$ via the relation $i(q, e) \neq 0$. Let $q_1, q_2,
q_3$ be three quadrilaterals in $\sigma$ so that $q_1 \to q_2 \to
q_3 \to q_1$ is the cyclic order induced by the cyclic order on
the opposite edges from a vertex. Let $W$ be the vector space with
a basis $\{q_1, q_2, q_3\}$. An anti-symmetric bilinear form
$\omega: W \times  W \to \bold R$ is defined by $\omega( q_i, q_j)
=1$ if and only if $(i, j) =(1,2),(2,3), (3,1)$. In particular,
$\omega(q_i, q_j) =-\omega(q_j, q_i)$. Given any two
quadrilaterals $q, q' \in \Box$, set $\omega(q, q')$ to be the
value just defined if they are in the same tetrahedron and
$\omega(q, q')=0$ if they are in different tetrahedra. In this
way, one obtains the Neumann-Zagier anti-symmetric bilinear form
$$ \omega: \bold R({\Box}) \times \bold R({\Box}) \to \bold R$$
where $\bold R(\Box)$ is the vector space with a basis $\Box$.
More details of the form can be found in the work of \cite{NZ},
\cite{Cho} and \cite{Til}. See also \cite{Lu2}.

The following was proved in \cite{NZ}.

\medskip
\noindent {\bf Proposition 3.1}(Neumann-Zagier). \it Suppose $(M,
\bold T)$ is a triangulated, oriented closed pseudo 3-manifold.
Then

(a) for any $q' \in \Box$, $\sum_{q \in \Box} \omega(q, q') =0$.

(b) for any pair of edges $e, e' \in E$,
$$ \sum_{ q, q' \in \Box}  i(q,e) i(q', e') \omega(q, q') =0.$$
\rm

\medskip
Indeed, part (a) follows from the anti-symmetric property, i.e.,
for any $i=1,2$ or $3$, $\sum_{j=1}^3  \omega(q_j, q_i) =0$. Part
(b) is more complicated. First, anti-symmetry shows that the
identity (b) holds if (1) $e =e'$, or (2) $e$ and $e'$ do not lie
in a tetrahedron, or (3) $e, e'$ lie in a tetrahedron and are
opposite edges. Now if $ e \neq e'$ and $e, e'$ lie in a
tetrahedron $\sigma$ and are not opposite, then $e, e'$ lie a
triangular face and there is a second tetrahedron $\sigma'$
containing $e, e'$. In this case, due to the orientations on
$\sigma$ and $\sigma'$, we have \begin{equation}\sum_{q, q'
\subset \sigma} i(q,e) i(q', e') \omega(q, q') =  -\sum_{q, q'
\subset \sigma'} i(q,e) i(q', e') \omega(q, q'). \end{equation}
Thus part (b) follows. For more details of the proof, see
\cite{NZ}, page 316-320.

\medskip

It is known (\cite{NZ})  that the restriction of the
Neumann-Zagier 2-form to the subspace $\{ x = \sum_{q \in \Box}
a_q q \in \bold R(\Box) |$ for each $\sigma \in  T$, $\sum_{ q
\subset \sigma} a_q =0$\} becomes  non-degenerated. The
2-dimensional counter-part of the Neumann-Zagier Poisson structure
is the Thurston's anti-symmetric bilinear form on the space of
measured laminations.  It is very closely related to the
Weil-Petersson symplectic form (\cite{PP}, \cite{BS}) on the
Teichmuller spaces and plays a vital rule in Kontsevich's work
(\cite{Ko}) on Witten's conjecture and many other works. It is
expected that Neumann-Zagier Poisson structure will play an
equally important role in (2+1) TQFT.

\medskip
\subsection{Thurston's gluing equation}

\noindent {\bf Definition 3.1} Suppose $(M, \bold T)$ is an
oriented closed pseudo 3-manifold with a triangulation and $k \in
(\bold S^1)^E$. Thurston's equation (with curvature $k$) is
defined for $ z \in \bold C^{\Box}$ so that for each $e \in E$,
\begin{equation} \prod_{q \in \Box} z(q)^{i(q, e)} =\pm k(z),
\end{equation}
and if $q,q' \in \Box$ so that $\omega(q,q') =1$, then
\begin{equation} z(q')(1-z(q))=1. \end{equation}

\medskip
By (3.3) and the fact that $f(t)=\frac{1}{1-t}$ satisfies $t
f(t)f(f(t))= -1$, we have, for each tetrahedron $\sigma \in T$,
\begin{equation}
\prod_{q \in \Box, q \subset \sigma} z(q) =-1. \end{equation}

Note that we do not require that $Im(z(q))
>0$ which corresponds to the positively oriented ideal tetrahedron
(\cite{NZ}). The work of Yoshida \cite{Yo} (see also \cite{NZ} and
\cite{Til}) shows that each solution $z$ so that the
right-hand-side of (3.2) is 1  produces a representation of
$\pi_1(M -V)$ to $PSL(2, \bold C)$.

Note that equation (3.2) is equivalent to
\begin{equation} \prod_{q \in \square} z(q)^{2i(q,e)} = k(e)^2 \end{equation}
It is the solution to (3.3) and (3.5) which is addressed in
theorem 1.1.

\subsection{Volume of $\bold S^1$-valued angle structures}

Recall that the Lobachevsky function $\Lambda(x) =-\int_0^x \ln |2
\sin(u)| du$ is a  continuous periodic function of period $\pi$
defined on $\bold R$.  It is real analytic on $\bold R -\pi \bold
Z$ so that $\lim_{ t \to 0} \Lambda '(t) = +\infty$. For more
details, see Milnor \cite{Mil}.  Given $t =e^{ \sqrt{-1} a} \in
\bold S^1$, define $\lambda(t) = \Lambda( a)$. This is well
defined since $\Lambda(a)$ has $\pi$ as a period. Furthermore,
$\lambda: \bold S^1 \to \bold R$ is real analytic on the subset
$\bold S^1 -\{\pm 1\}$. For an $\bold S^1$-valued angle structure
$x : \square \to \bold S^1$ on $(M, \bold T)$, define its \it
volume \rm $\bold V(x)$ to be
$$ \bold V(x) = \sum_{ q \in \square} \lambda( x(q)) =\sum_{ q \in \Box} \Lambda(\arg(x(q))). $$
The volume function $\bold V$ is continuous and in particular, has
a maximum and a minimum point. By definition the smooth points for
$\bold V: SAS( \bold T, k) \to \bold R$ are exactly those points
$x$ where $x(q) \neq \pm 1$ for all $q$.

\medskip
The main theorem in the paper can be stated as,

\medskip
\noindent {\bf Theorem 3.2} \it Suppose $(M, \bold T)$ is
triangulated oriented closed pseudo 3-manifold  and $k \in (\bold
S^1)^E$ satisfying (2.16) and (2.17).
 Let $p$ be a maximum point
of the volume function $\bold V :SAS( \bold T, k)  \to \bold R$.

(a) If $p$ is a smooth point for $\bold V$, then $p$ produces a
solution to the generalized Thurston gluing equation (3.3) and
(3.5).

(b) If $p$ is a non-smooth point for $\bold V$, i.e., $p(q_0) =
\pm 1$ for some $q_0 \in \Box$, then $p$ produces a 2-quad-type
solution $y$ to Haken's normal surface equation so that $y(q_0)
\neq 0$.  \rm
\medskip

\subsection{Smooth critical point of the volume}

The following lemma was known to Casson and Rivin.

\medskip
\noindent {\bf Lemma 3.3} \it Suppose $x \in SAS(\bold T, k)$ is a
smooth critical point of the volume $\bold V: SAS(\bold T) \to
\bold R$.  Then  Thurston's equation (3.3) and (3.5) has a
solution in $(\bold C -\bold R)^{\Box}$. \rm
\medskip

\medskip
\noindent {\bf Proof.} Suppose $q_1, q_2, q_3$ are three
quadrilaterals in a tetrahedron.  Let $x_i=x(q_i)$ be the $\bold
S^1$-valued angle at the quadrilateral. We define the associated
complex values $z(q_i)$ by the formula,
$$ z(q_i) = \frac{ x_j - \bar x_j}{x_k - \bar x_k}x_i =\frac{ \sin(\arg(x_j))}{\sin(\arg(x_k))} x_i
$$ where $\omega(q_i, q_j) =1$ and $\{i,j,k\}=\{1,2,3\}$.
This is well defined since $x_k-\bar x_k \neq 0$ by the
assumption. More generally, for $x \in SAS(\bold T)$ and $x(q)
\neq \pm 1$ for all $q$, one defines $z \in \bold C^{\Box}$, by
$$ z(q)  = x(q) \prod_{ r \in \Box}  (\sin(\arg(x(r))))^{
\omega(r,q)}.$$

We claim that $z$ is a solution to Thurston's equation (3.3) and
(3.5).

First, (3.3) follows by  a direct calculation and the definition.
Let us assume that $z_i = z(q_i)$  and that $\omega(q_1, q_2) =1$.
By definition, we have
$$ z_1 = \frac{ x_2-\bar x_2}{x_3-\bar x_3} x_1$$
and $$ z_2 = \frac{ x_3-\bar x_3}{x_1-\bar x_1} x_2.$$ Due to
$x_1x_2x_3=-1$ and $x_i \bar x_i =1$,  then (3.3) says that
$$ z_2(1-z_1)=1.$$
Indeed,
$$z_2( 1-z_1) = (\frac{x_3-\bar x_3}{x_1 -\bar x_1} x_2) \frac{ x_3 -\bar x_3 - x_1 x_2 + x_1 \bar x_2}{ x_3
-\bar x_3} = \frac{ x_3 +x_1 \bar x_2}{ x_1 -\bar
x_1}x_2=\frac{x_3 x_2 +x_1}{x_1 -\bar x_1} = 1$$

To see (3.5), we need to use the critical point equation for
$\bold V$ at the smooth point $x$. By definition, we can identify
the tangent space to a point of $SAS(\bold T, k)$ with $TAS(\bold
T)$. Indeed, for any $v \in TAS(\bold T)$ and $ x \in SAS(\bold T,
k)$, the path $ p(t)=x e^{tv} \in SAS(\bold T, k)$ given by
$$ p(t)(q) =  x(q) e^{ t v(q)}$$
for $ t \in (-\epsilon, \epsilon)$ has tangent vector $v$ at $t=0$
and all tangent vectors to $SAS(\bold T, k)$ at $x$ are of this
form. Now due to $x(q) \neq \pm 1$, $\frac{d \bold V( x
e^{tv})}{dt}|_{t=0} =0$ shows that,
\begin{equation}
\sum_{ q \in \Box}  v(q) \ln|\sin ( \arg(x(q)))| =0.
\end{equation}

Choose a specific $v \in TAS(\bold T)$ as follows. Fix an edge $e
\in E$, by proposition 3.1,  \begin{equation} v_e = \sum_{q \in
\Box} \sum_{r \in \Box} i(q, e) \omega(r, q)  r^* \in TAS(\bold
T). \end{equation}  Now substitute $v_e$ for $v$ in (3.6) and use
the fact that $\sum_{ q \in \Box} i(q,e) v_e(q) =0$, we obtain,
for each $e \in E$,
$$  \prod_{ q \in \Box} z(q)^{i(e,q)}
= \prod_{ q \in \Box} x(q)^{ i(e, q)}   \prod_{r \in \Box}
\sin(\arg(x(q)))^{ i(e,q) \omega(r, q)}$$
$$ =  k(e) \prod_{ q, r \in \Box} \sin(\arg(x(q)))^{ i(e,q) \omega(r, q)}=\pm k(e) $$
due to (3.6) and (3.7). This verifies (3.5) and ends the proof.

\medskip

Furthermore, if $z$ is a solution in $(\bold C -\bold R)^{\Box}$
to generalized Thurston's equation (3.3) and (3.5) over a closed
3-manifold, then $\frac{z}{|z|}$  is an $\bold S^1$-angle
structure which is a smooth critical point of the volume $\bold V$
in $SAS(\bold T)$. The proof uses the fact that for a closed
manifold $M$, the tangent space $TAS(\bold T)$ is generated by the
vectors $v_e$'s for $e \in E$ given by (3.7). We omit the details.

\medskip

\noindent {\bf Corollary 3.4} \it Under the same assumption as in
lemma 3.3, if $x \in SAS(\bold T, k)$ is a smooth critical point
of the volume $\bold V$, let $y \in \bold R_{\geq 0}^{\Box}$ be
the vector so that $y(q) = - \ln |(\sin( \arg(x(q)))|$. Then $y
\in Proj_{\Box}(\bold S_{ns})$. \rm

\medskip
Indeed, (3.6) shows that $ y \in TAS(\bold T)^{\perp}$. Thus, by
theorem 2.2, $ y \in Proj_{\Box}(\bold S_{ns})$.

It will be very interesting to see what the topological
information $y$ contains.

\section{Volume optimization and normal surfaces}

A relationship between the smooth critical points of the volume
$\bold V: SAS(\bold T,k) \to \bold R$ and the normal surfaces is
established in Corollary 3.4. In this section, we will investigate
the case of non-smooth critical points of the volume.
 Since the
function $\bold V$ is not smooth, the definition of the critical
points of $\bold V$ should be specified.  First of all, we will
show (corollary 4.3) that for any $p \in SAS(\bold T)$ and $u \in
TAS(\bold T)$ the  limit $\lim_{t \to 0} \frac{d \bold V(p
e^{tu})}{dt}$ always exists as an element in $[-\infty, \infty] =
\bold R \cup \{ \infty, -\infty \}$. We say that a point $ p \in
SAS(\bold T,k)$ is a \it critical point \rm of the volume $\bold
V$ if for all $u$ in $TAS(\bold T)$,
\begin{equation}
\lim_{t \to 0}  \frac{d \bold V(p e^{tu})}{dt} =0.
\end{equation}

Using this definition, one sees easily that the maximum and
minimum points of $\bold V$ are critical points, i.e., the volume
function $\bold V$ always has critical points.

The main theorem in the section, which implies theorem 1.1, is the
following,

\medskip
\noindent {\bf Theorem 4.1} \it Suppose $(M, \bold T)$ is an
orientable closed triangulated pseudo 3-manifold with $SAS(M,
\bold T,k) \neq \emptyset$. If the volume $\bold V: SAS(M, \bold
T, k) \to \bold R$ has a non-smooth critical point $p$, i.e.,
$p(q_0) = \pm 1$ for some $q_0 \in \Box$, then $p$ produces a
2-quad-type solution $y$ to Haken's normal surface equation so
that $y(q_0) \neq 0$.  \rm
\medskip

Recall that by proposition 2.6, $SAS(M, \bold T,k) \neq \emptyset$
if and only if $k$ satisfies (2.16) and (2.17).

\subsection{ Subderivatives of the volume function}

The volume function $\bold V$ is essentially composed by the
function $W: P \to \bold R$ where $P=\{(x_1, x_2, x_3) \in \bold
R^3 | x_1 + x_2 + x_3 =\pi\}$ and $W(x_1, x_2, x_3) =
\Lambda(x_1)+\Lambda(x_2) + \Lambda(x_3)$.  The function $W$ is
not smooth on the subset defined by some $x_i \in \pi \bold Z$.
However, we can obtain subderivative information of $W$ at these
points.

The function $h(t) = t \ln |t|$ can be extended to be a continuous
function from $\bold R \to \bold R$ by declaring $h(0) =0$. In the
sequel, this extension, still denoted by $t \ln|t|$, will be used.

\medskip
\noindent {\bf Lemma 4.2} \it Take a point $a=(a_1, a_2, a_3) \in
P$ and $b =(b_1, b_2, b_3) \in \bold R^3$ so that $b_1+b_2+b_3=0$.
Define $f(t) = \frac{ d W(a+tb)}{dt}$.  Then $\lim_{t \to 0} f(t)$
exists as an element in $\bold R \cup \{\pm \infty\}$ and

(a) if $a_i \notin \pi \bold Z$ for all $i$,
\begin{equation}
\lim_{t \to 0} f(t) =-\sum_{i=1}^3 b_i \ln |\sin(a_i)|,
\end{equation}

(b) if $a_i \in \pi \bold Z$ for all $i$, \begin{equation} \lim_{t
\to 0} f(t) = -\sum_{i=1}^3 b_i \ln | b_i|
\end{equation}

(c) if $a_1 \in \pi \bold Z$ and $a_2, a_3 \notin \pi \bold Z$,
then
\begin{equation}
\lim_{t \to 0} (f(t) + b_1 \ln|t|)  = -b_1 \ln|b_1| -\sum_{ i=2}^3
b_i \ln | \sin (a_i)|
\end{equation}
\rm
\medskip

\noindent {\bf Proof.} We have $f(t) = -\sum_{i=1}^3 b_i \ln|2
\sin(a_i+t b_i)| = -\sum_{i=1}^3 b_i \ln|\sin(a_i+ tb_i)|$ due to
$\sum_{i=1}^3 b_3=0$. Now part (a) follows form the definition.

For part (b), due to $\ln(| \sin( t + \pi)|) = \ln|\sin(t)|$, it
follows that $f(t) = -\sum_{i=1}^3 b_i \ln(|\sin(tb_i)|)$. The
result is obvious if $b_i=0$ for all $i$. Otherwise, say $b_3 \neq
0$, then $b_3 =-b_1-b_2$.  Substitute it to $f(t)$, we obtain
 $$ f(t)
= -b_1 \ln|\frac{\sin(b_1t)}{\sin(b_3t)}|-b_2 \ln|\frac{\sin(b_2
t)}{\sin(b_3 t)}|.$$  By taking the limit as $t \to 0$, we obtain
part (b).

For part (c), we write $$f(t) = -b_1 \ln| \frac{\sin(b_1
t)}{b_1t}| -b_1 \ln|b_1 t| -\sum_{i=2}^3 b_i \ln | \sin( a_i+ t
b_i)|$$
$$= -b_1 \ln|t| -b_1 \ln |b_1| -\sum_{i=2}^3 b_i \ln| \sin(a_i)| +
o(t).$$ where $o(t)$ is the quality so that $\lim_{t \to 0} o(t)
=0$. This establishes part (c) and finishes the proof.

\medskip
Not that due to $a_1+a_2+a_3=\pi$, cases (a), (b) and (c) are the
list of all cases up to symmetry.
 The limit $\lim_{t \to 0} f(t)$
in cases (b), (c) above is called the subderivative of the
function $W$ at the point $a$. The subderivative, considered as a
function of the tangent vector $b$, is homogeneous of degree 1.
However, due to the non-smoothness, the subderivative, as shown in
(b), (c), is not a linear function of $b$.

In the case of the $\bold S^1$-valued angle structure, consider $X
=\{ a=(a_1, a_2, a_3) \in (\bold S^1)^3 | a_1a_2a_3 =-1\}$ and the
volume $\bold V(a) = \sum_{ i=1}^3 \lambda(a_i) =\sum_{i=1}^3
\Lambda(\arg(a_i))$. Consider a tangent vector  $b =(b_1, b_2,
b_3) \in \bold R^3$ so that $b_1+b_2+b_3=0$. Define $f(t) = \frac{
d \bold V(ae^{tb})}{dt}$. Then $\lim_{t \to 0} f(t)$ exists as an
element in $\bold R \cup \{\pm \infty\}$ and  the above lemma
says,

(a) if $a_i \neq  \pm 1$ for all $i$,
\begin{equation}
\lim_{t \to 0} f(t) =-\sum_{i=1}^3 b_i \ln |\sin( \arg(a_i))|,
\end{equation}

(b) if $a_i = \pm 1$ for all $i$, \begin{equation} \lim_{t \to 0}
f(t) = -\sum_{i=1}^3 b_i \ln | b_i|
\end{equation}

(c) if $a_1 =\pm 1$ and $a_2, a_3 \neq \pm 1$, then
\begin{equation}
\lim_{t \to 0} (f(t) + b_1 \ln|t|)  = -b_1 \ln|b_1| -\sum_{ i=2}^3
b_i \ln | \sin (\arg(a_i))|
\end{equation}

\medskip
\noindent{\bf Corollary 4.3} \it For any $a \in SAS(\bold T, k)$,
 there exists a linear function $g(b)$ of
$b \in TAS(\bold T)$ and a continuous function $f(b,t)$ of $b$ and
$t \in (-\epsilon, \epsilon)$ so that
$$ \frac{ d \bold V( a e^{ tb})}{ dt} = g(b) \ln |t| + f(b,t).$$
In particular, the limit $\lim_{ t \to 0} \frac{ d \bold V( a e^{
tb})}{ dt}$ always exists as an element in $[-\infty, \infty]$.
Furthermore,  a local maximum or minimum point of $\bold V$ is a
critical point. \rm

\subsection{ A proof of theorem 4.1}

Suppose $ p \in SAS(\bold T,k)$ is a non-smooth critical point of
the volume function $\bold V$ so that $p(q_0) =\pm 1$ for some
$q_0 \in \Box$. By definition of critical points, we have $$\lim_{
t \to 0} \frac{ d \bold V( p e^{tb})}{ dt} =0$$ for all $b$ in
$TAS(\bold T)$. By the definition of $\bold V$, we have
$$ \bold V(x) =\sum_{ \sigma \in T} \sum_{ q \in \Box, q \subset \sigma}
\Lambda( arg(x(q))).
$$
Let $Y =\{ q \in \Box | p(q) =\pm 1 \}$ which contains $q_0$ and
$Y' =\{ q \in Y | $ there exists $\sigma \in T$ so that $q \subset
\sigma$ and for the other two $q', q'' \subset \sigma$ $p(q'),
p(q'') \neq \pm 1$\}.

Thus, by (4.5)-(4.7), we can write
$$
\lim_{t \to 0} (\frac{ d\bold V(p e^{tb})}{dt} + \ln |t| \sum_{q
\in Y'} b(q)) = -\sum_{ q \in Y} b(q) \ln|b(q)|  -\sum_{ q \notin
Y} b(q) \ln|\sin(\arg(p(q)))|. $$

By the critical point condition (4.1), we obtain $\sum_{q \in Y'}
b(q)=0$ and
\begin{equation}
\sum_{ q \in Y} b(q) \ln|b(q)|= -\sum_{q \notin Y} b(q) \ln |
\sin(\arg(p(q)))|
\end{equation}
For each $q \in \Box$, let $f_q: TAS(\bold T) \to \bold R$ be the
linear function on $TAS(\bold T)$
 defined by $f_q(b) = b(q)$.   Then the right-hand-side of (4.8) is a linear function in $b$
on $TAS(\bold T)$ and the left-hand-side of (4.8) is a sum of the
functions $ f_q(b) \ln |f_q(b)|$.

\medskip\noindent {\bf Lemma 4.4} \it Suppose $W$ is a finite
dimensional vector space over $\bold R$ and $f_1, ..., f_n$ and
$g$ are linear functions on $W$ satisfying \begin{equation} \sum_{
i=1}^n f_i(x) \ln | f_i( x)| = g(x). \end{equation} Then for each
index $i$ there exists $ j \neq i$ and $\lambda_{ij} \in \bold R$
so that
$$f_i(x) = \lambda_{ij} f_j(x).$$ \rm

\medskip
\noindent {\bf Proof.}  Note if one of $f_i$ is the zero function
$f_i(x) =0$ for all $x$, then the lemma holds. Let us assume that
all $f_i$'s are non-zero functions. We may assume that $W = \bold
R^m$ and $x=(x_1, ..., x_m) \in W$ after a linear change of
variables. Write
$$ f_i(x) = \sum_{j=1}^m a_{ij} x_j.$$ Now suppose the result does
not hold, say $f_1(x)$ is not propositional to $f_j(x)'s$ for $j
\geq 2$. Then we can find a point $v \in ker(f_1)$ so that $ v
\notin \cup_{j=2}^n ker(f_j)$. Since $f_1 \neq 0$, for simplicity,
let us assume that $a_{11} \neq 0$. Now take derivative of (4.9)
with respect to $x_1$. We obtain an equation of the form
\begin{equation}
\sum_{j=1} a_{1j} \ln| f_j(x)| =  h(x)
\end{equation} where $h(x)$ is  a linear function.  Take a
sequence of vectors $x$ converging to $v$ in (4.10), we obtain a
contradiction since $a_{11} \neq 0$.  This ends the proof.

\medskip

Applying lemma 4.4 to (4.8) with $f_i$'s being $f_q$'s for $q \in
Y$, we conclude that for $f_{q_0}$, there exist $f_{q_1}$, $q_1
\neq q_0$ and $\lambda \in \bold R$ so that $f_{q_0}(b) = \lambda
f_{q_1}(b)$ for all $b \in TAS(\bold T)$.  This shows that $(b,
q_0^* - \lambda q_1^*) =0$ for all $b$. By theorem 2.2, there
exists a solution $y$ to Haken's equation so that $y(q_0)=1,
y(q_1) = -\lambda$ and $y(q) =0$ for all $ q \in \Box-\{ q_0,
q_1\}$. This ends the proof.

\medskip
\subsection{A generalization}

In an unpublished work \cite{FG}, David Futer and Fran{\c c}ois
Gu\'eritaud proved a very nice theorem concerning the non-smooth
maximum points of the volume function. Given $x \in SAS(\bold T)$,
we say a tetrahedron $\sigma \in T$ is \it flat \rm with respect
to $x$ if $x(q) = \pm 1$ for all $q \subset \sigma$ and \it
partially flat \rm if $x(q) = \pm 1$ for one $q \subset \sigma$.

\medskip
\noindent {\bf Theorem 4.5} (Futer-Gu\'eritaud) \it Suppose $(M,
\bold T)$ is an oriented triangulated closed pseudo 3-manifold. If
$x$ is a non-smooth maximum point of the volume function on
$SAS(\bold T)$, then there exists a non-smooth maximum volume
point $y \in SAS(\bold T)$ so that all partially flat tetrahedra
in $y$ are flat. \rm

\medskip
A written proof it, supplied by Futer-Gu\'eritaud, can be found in
\cite{Lu2}. Combining theorems 4.5 and 4.1, we obtain a stronger
statement that

\medskip
\noindent {\bf Theorem 4.6} \it Suppose $(M, \bold T)$ is a closed
triangulated oriented pseudo 3-manifold so that it has a
non-smooth maximum volume point in $SAS(\bold T)$. Then there
exist three 2-quad type solutions $x_1, x_2, x_3$ of Haken's
normal surface equation so that there are three distinct
quadrilaterals $q_1, q_2, q_3$ in a tetrahedron with $x_i(q_i)
\neq 0$ for all $i$. \rm

\medskip

We call three 2-quad-type solutions appeared in theorem 4.6 a \it
cluster of 2-quad-type solutions. \rm  In the joint work
\cite{LT2}, Tillmann and I proved the following topological
result.

\medskip
\noindent {\bf Theorem 4.7} (\cite{LT2}) \it
 Suppose $(M, T)$ is a minimally triangulated orientable closed
3-manifold which supports a cluster of three 2-quad-type solutions
to Haken's equation. Then,

(a) $M$ is reducible, or

(b) $M$ is toroidal, or

(c) $M$ is a Seifert fibered space, or

(d) $M$ contains the connected sum $\#_{i=1}^3 RP^2$ of three
copies of the projective plane. \rm

\medskip
Theorem 1.2 mentioned in the introduction is a consequence of
theorems 4.6 and 4.7.

\medskip
\section{Some questions}

Based on theorem 4.7, we propose the following conjecture which is
weaker than conjecture 1.

\medskip
\noindent {\bf Conjecture 2} Suppose $(M, \bold T)$ is a minimally
triangulated irreducible orientable closed 3-manifold so that the
triangulation does not have a cluster of three 2-quad-type
solutions to Haken's equation. Then there is a solution to
Thurston equation associated to $\bold T$.

\medskip
Note that by the same argument as in \S1.2 and using theorem 4.7
instead of theorem 1.2 and Segerman-Tillmann's theorem, we see
that conjecture 2 for simply connected manifold is equivalent to
the Poincar\'e conjecture.

Conjecture 2 relates solutions of Thurston equation to that of
Haken equation and does not involve volume optimization process.
The minimality condition in conjecture 2 is necessary. This was
shown to us by Ben Burton and Henry Segerman.  Solutions to
Thurston's equation have been found in many cases. For instance,
Tillmann proved in \cite{Til2} that if $M$ is a non-compact finite
volume hyperbolic 3-manifold and $\bold T$ is an ideal
triangulation so that each edge is homotopically essential, then
$\bold T$ supports a solution to Thurston's equation. However, a
general existence theorem for solving Thurston's equation seems to
be still lacking. Conjecture 2 is an attempt to address the issue.

\medskip
Given a solution to Thurston equation, by the work of Yoshida
\cite{Yo}, one can produce a representation of the fundamental
group of $M -V$ to $PSL(2, \bold C)$ where $V$ is the set of
vertices. It is  interesting to know when solutions to Thurston's
equation produce irreducible representations of the fundamental
group. See the work of \cite{Fra} and \cite{FK}.

Finally solving Thurston's equation over the real numbers, i.e.,
$z \in \bold R^{\Box}$, seems to be an attractive problem.  Here
is a step toward producing a real-valued solution to Thurston
equation.

\medskip
\noindent {\bf Definition 5.1}. Let $\bold Z_2$ be the field of
two elements $\{0,1\}$. A \it $\bold Z_2$-taut structure \rm on a
triangulated closed pseudo 3-manifold $(M, \bold T)$ is a map $f:
\Box \to \{0,1\}$ so that

(a) if $q_1, q_2, q_3$ are three quadrilaterals in each
tetrahedron $\sigma$, then $\{f(q_1), f(q_2), f(q_3)\}
=\{0,0,1\}$, and

(b) for each edge $e$ in $\bold T$, $\sum_{ q \in \Box} i(q,e)
f(q) =0$.

\medskip

The motivation for the definition comes from taut triangulations
and real-valued solutions to Thurston's equation. Indeed, if $z$
is a real-valued solution to Thurston's equation, then there is an
associated $\bold Z_2$-taut structure $f$ defined by:  $f(q) =0$
if $z(q) >0$ and $f(q)=1$ if $z(q) < 0$.  Another motivation comes
from taut triangulation. Suppose $\bold T$ is a taut
triangulation, i.e., there is a map $g: \Box \to \{0, \pi\}$ so
that for each tetrahedron $\sigma$, $\sum_{ q \subset \sigma} g(q)
= \pi$ and for each edge $e$, $\sum_{q} i(q,e) g(q) = 2\pi$. Then
one defines a $\bold Z_2$-taut structure by $f(q) = \frac{1}{\pi}
g(q)$.  A very interesting question is to find condition on $\bold
T$ so that $\bold Z_2$-taut structures exist. Is it possible that
the non-existence of $\bold Z_2$-taut structures implies the
existence of some special solutions to Haken's normal surface
equation?

Finally, Tillmann and I observed that the equations for $\bold
Z_2$-taut structures are non-linear but quadratic in $f(q)$.
Indeed, a vector $f \in \bold Z_2^{\Box}$ is a $\bold Z_2$-taut
structure if and only if condition (b) in definition 5.1 holds and
for each tetrahedron $\sigma$
\begin{equation}
 \sum_{ q \subset \sigma} f(q) =1, \end{equation} and
$$\sum_{q \neq q', q,q' \subset \sigma} f(q) f(q') =0.$$

The condition (b) in definition 5.1 and (5.1) should be considered
as the definition of a $\bold Z_2$-angle structure.

\medskip
\section{Appendix}

We give a new proof of Kang-Rubinstein theorem in this section.
First, one checks easily that both $W_{\sigma}$ and $W_e$ are in
$\bold S_{ns}$. Next, by a simple dimension counting, one sees
that $\dim (\bold S_{ns}) \leq |E| + |T|$. Thus, it suffices to
prove that $\{ W_{\sigma}, W_e | \sigma \in T, e \in E\}$ is an
independent set. To this end, suppose otherwise that there exists
$ h \in \bold R^E \times \bold R^T$ so that
$$ \sum_{e \in E} h(e) W_e + \sum_{ \sigma \in T} h(\sigma)
W_{\sigma} =0.$$

We can write it as, $$ \sum_{ t \in \triangle} (-\sum_{ e > t}
h(e)- \sum_{\sigma > t} h(\sigma)) t^*  + \sum_{q \in \square}
(\sum_{ e \in E} h(e) i(q, e) +\sum_{ \sigma \in T, q \subset
\sigma} h(\sigma)) q^*  =0 $$

Since $\{t^*,q^*\}$ form a basis, we obtain for each $t \in
\triangle$
\begin{equation}
\sum_{ e > t} h(e)- \sum_{\sigma > t} h(\sigma) =0
\end{equation}
and for each $q \in \square$
\begin{equation}
\sum_{ e \in E} h(e) i(q, e) +\sum_{ \sigma \in T, q \subset
\sigma} h(\sigma)=0
\end{equation}

Consider a fixed tetrahedron $\sigma \in T$. We claim that the
system of linear equations (6.1) and (6.2) for the six edges of
$\sigma$ has only the trivial solution, i.e., $h(e) = h(\sigma)
=0$. In particular, this shows that $\{W_e, W_{\sigma} \}$ is
independent.

To see the claim, let us label the vertices of $\sigma$ by
$1,2,3,4$ and the six edges by $e_{ij}$ were $ i \neq j \in
\{1,2,3,4\}$. Let $h_{ij} = h(e_{ij})$ and $f=h(\sigma)$. Then
(6.1) and (6.2) say: at the i-th vertex \begin{equation}
h_{ij}+h_{ik}+h_{il}=f \end{equation} and
\begin{equation} h_{ij}+h_{kl}= f \end{equation} for $\{i,j,k,l\}=\{1,2,3,4\}$.
Consider the sum of two equations (6.3) at the i-th and j-th
vertices subtracting the sum of the two equations (6.3) at the
k-th and l-th vertices. We obtain, $h_{ij}=h_{kl}$, i.e.,
$h(e)=h(e')$ when $e,e'$ are opposite edges. Now by (6.4), we see
that $h_{ij} = f/2$ for all $i \neq j$. Now substitute back to
(6.3), we obtain $3f/2 = f$. Thus $f=0$ and $h_{ij}=0$, i.e.,
$h(e)=h(\sigma)=0$.

Department of Mathematics

Rutgers University

Piscataway, NJ 08854

email: fluo\@math.rutgers.edu
\end{document}